\documentclass[11pt]{article}

\usepackage{lscape,tikz}
\usepackage{authblk}
\definecolor{darkblue}{rgb}{0.2,0.2,0.71}
\usepackage{todonotes}
\usepackage{epsfig,cancel,ulem, amsthm,amssymb}
\usepackage{color,soul}
\usetikzlibrary{decorations.markings,arrows}

\usepackage{wrapfig}
\usepackage{color}
\usepackage{graphicx,framed}
\usepackage[colorlinks=true, pdfstartview=FitV, linkcolor=darkblue, citecolor=darkblue, urlcolor=darkblue]{hyperref}
\definecolor{shadecolor}{rgb}{0.95, 0.95, 0.86}
\definecolor{darkgreen}{rgb}{0.2, 0.5,  0}

\newcommand{\rk}{\mathrm{rank}}

\def\&{\vspace{-5pt}&}

\def \eqref#1{(\ref{#1})}
\def \& {&\hspace{-10pt}}

\newcommand{\bt}{\beta}

\renewcommand{\O}{\Omega}

\newcommand{\br}{{\mathbb R}}

\newcommand{\g}{\mathfrak{g}}

     \newcommand{\lop}[1]{\mathfrak{L}(#1)}
      \newcommand{\Lie}{\mathfrak{L}}
\newcommand{\bil}[2]{{\langle #1 | #2\rangle}}

     \newcommand{\W}{\mathcal{W}}

\newcommand{\f}{\mathcal{F}}
\newcommand{\I}{\mathcal{I}}
\newcommand{\Q}{\mathcal{Q}}
\newcommand{\N}{\mathcal{N}}
 
\newtheorem{theorem}{Theorem}[section]
\newtheorem{example}[theorem]{Example}
\newtheorem{exercise}[theorem]{Exercise}

\newtheorem{lemma}[theorem]{Lemma}
\newtheorem{remark}[theorem]{Remark}

\newtheorem{proposition}[theorem]{Proposition} 
\newtheorem{corollary}[theorem]{Corollary} 
\newtheorem{definition}[theorem]{Definition}

\def\bt{\begin{theorem}}
\def\et{\end{theorem}}
\def\bc{\begin{corollary}}
\def\ec{\end{corollary}}
\def\bx{\begin{example}\small}
\def\ex{\end{example}}
\def\bxr{\begin{exercise}\small}
\def\exr{\end{exercise}}
\def\bl{\begin{lemma}}
\def\el{\end{lemma}}
\def\bd{\begin{definition}}
\def\ed{\end{definition}}
\def\bp{\begin{proposition}}
\def\ep{\end{proposition}}

\def\br{\begin{remark}}
\def\er{\end{remark}}

\def\be{\begin{equation}}
\def\ee{\end{equation}}

\def\beq{\begin{equation}}
\def\eeq{\end{equation}}
\def\&{\hspace{-15pt}&}
\def\bea{\begin{eqnarray}}
\def\eea{\end{eqnarray}}

\def\1{{\bf 1}}

\newcommand{\h}{\mathfrak{h}}

\newcommand{\ad}{\mathrm{ad}}

\makeatletter
\@addtoreset{equation}{section}
\makeatother

\begin{document}

\title{Low dimensional bihamiltonian structures of topological type}
\author{Yassir Ibrahim Dinar}
\affil{Sultan Qaboos University, Oman}
\date{}
\maketitle




\begin{abstract}
 We construct  local bihamiltonian structures  from  classical  $W$-algebras  associated to  non-regular nilpotent elements of regular semisimple type in Lie algebras of type $A_2$ and $A_3$.  They form  exact Poisson pencil, admit a dispersionless limit and their    leading terms define logarithmic or trivial Dubrovin-Frobenius manifolds. We  calculate the  corresponding central invariants which are expected to be constants. In particular, we get  Dubrovin- Frobenius manifolds associated to the focused Schr\"{o}dinger equation and   Hurwitz space $M_{0;1,0}$ and the corresponding bihamiltonian  structures of topological type.   
\end{abstract}

\maketitle

\vspace{1.5cm}
\section{Introduction}

Dubrovin-Frobenius manifold is a marvelous geometric realization introduced by Boris Dubrovin for a potential satisfying an undetermined  partial differential equations known as Witten-Dijkgraaf-Verlinde-Verlinde (WDVV) equations which describe the  module space of two dimensional topological field theory \cite{DuRev}. One of the  main methods to construct   Dubrovin-Frobenius manifolds  exists within the theory of  local bihamiltonian structure of hydrodynamic type, i.e., local compatible Poisson structures of hydrodynamic type \cite{DFP}. Moreover,  one can obtain a local bihamiltonian structure of hydrodynamic type by taking, under certain truncation, the leading term of a  local bihamiltonian structure admitting a dispersionless limit. On the other hand, if the leading term of a bihamiltonian structure  defines a semisimple Dubrovin-Frobenius manifold, then one can calculate the central invariants of the bihamiltonian structure. Besides, if the central invariants are constant and equal then the bihamiltonian structure is of topological type, i.e., theoretically one can construct it  from  the Dubrovin-Frobenius manifold structure using certain identities  inspired by theory of Gromov-Witten invariants \cite{DZ}. 

In recent work \cite{mypaper6}, we gave a uniform construction of algebraic Dubrovin-Frobenius manifolds beginning from Drinfeld-Sokolov  bihamiltonian structures associated to certain type of nilpotent elements in simple Lie algebras. Recall that one of the compatible local Poisson brackets is known as (or satisfies identities leading to) classical $W$-algebra \cite{fehercomp}.   More precisely, we fix a distinguished nilpotent element of regular semisimple type in a simple Lie algebra  of rank 
$r$. By definition, distinguished means its adjoint orbit  has no representative in a proper Levi subalgebra \cite{COLMC} while  of regular semisimple type indicates that one can construct a certain  Cartan subalgebra  known as the opposite Cartan subalgebra \cite{Elash}. Then, we considered  the space of common equilibrium points $N$ defined by the leading term of the associated  Drinfeld-Sokolov  bihamiltonian structure. The space $N$ is of dimension $r$ and we used  Dirac reduction to get on its loop space a new local bihamiltonian structure. This  new bihamiltonian structure  admits a dispersionless limit and its leading term gives  an algebraic Dubrovin-Frobenius manifold. These uniform construction   depends  on the associated opposite Cartan subalgebra.  Then it is natural to ask about the possibility to generalize the  construction  \cite{mypaper6} for arbitrary nilpotent elements and to inquire about the types of resulting Dubrovin-Frobenius manifolds.

 In this article, we consider   non-distinguished  nilpotent elements  of regular semisimple type  in simple Lie algebras $sl_3$ and $sl_4$ (they form 3 types of nilpotent elements) to begin answering the above-mentioned questions. For any nilpotent element under consideration,  the associated Drinfeld-Sokolov bihamiltonian structure is reduced to   a bihamiltonian structure admitting a dispersionless limit  on the loop space of the space of common equilibrium points but  its leading term fails to define a Dubrovin-Frobenius manifold structure using the methods given in \cite{mypaper6}.  However, starting from the classical $W$-algebra, we find  another bihamiltonian structure  which forms an exact Poisson pencil and  admits a dispersionless limit and its leading term defines a logarithmic, algebraic or a trivial Dubrovin-Frobenius manifold. We also calculate the central invariants of the new bihamiltonian structures which are expected to be constant (it follows from the exactness of  the bihamiltonian structure \cite{FalLor}). For a subregular nilpotent element, the  central invariants of the bihamiltonian structure  are equal. Moreover, in the case of $sl_3$ (resp. $sl_4$), the Dubrovin-Frobenius manifold  structure is associated to the focused Schr\"{o}dinger equation \cite{DF} (resp. is constructed on Hurwitz space $M_{0;1,0}$ of  genius 0 with two marked points of degree 2 and 1 \cite{DuRev}). Thus, we obtain  from classical $W$-algebra the bihamiltonian structure of topological type for logarithmic Dubrovin-Frobenius manifolds.   

The paper is organized as follows. In section 2 and 3, we recall the notion of Dubrovin-Frobenius manifolds and Drinfeld-Sokolov reduction. Then in  section \ref{regular type},  we  review  the construction of algebraic Dubrovin-Frobenius manifolds given in  \cite{mypaper6}. In the last section, we consider in details the case of  non-distinguished nilpotent elements of regular semisimple type in Lie algebra $sl_3$ and $sl_4$. 

\section{Dubrovin-Frobenius manifolds}

Let $M$ be a smooth manifold of dimension $r$ and  fix  local coordinates $(u^1, . . . , u^r)$ on $M$. In what follows, summation with respect to repeated upper and lower indices is assumed. From \cite{DuRev} and \cite{DFP}, we recall the  following definitions.

\bd \label{contra metric}  A symmetric bilinear form $(. ,. )$ on $T^*M$ is called a contravariant
metric if it is invertible on an open dense subset $M_0 \subseteq M$. We define the contravariant Christoffel symbols $\Gamma^{ij}_k$  for a contravariant
metric $(. ,. )$ by
\[
\Gamma^{ij}_k:=-\Omega^{im} \Gamma_{mk}^j
\]
where $\Gamma_{mk}^j$ are the  Christoffel symbols of the metric $<. ,. >$ defined on $TM_0$ by the inverse of the matrix $\Omega^{ij}(u)=(du^i, du^j)$.
We say the metric $(.,.)$ is flat if  $<. ,. >$ is flat.
\ed

Note that the Christoffel symbols given in definition \ref{contra metric} determine for the metric  $\O$ the contravariant (resp. covariant) derivative  $\nabla^{i}$ (resp.  $\nabla_{i}$) along the covector $du^i$ (resp. the vector field $\partial_{u^i}$). They are related by the identity  $\nabla^{i}=\O^{ij}(u) \nabla_{j}$. We denote the Lie derivative of $\Omega$ along a vector field $V$ by $\Lie_V\Omega$.

\bd 
A flat pencil of metrics (FPM)  on $M$ is a pair $(\Omega_2,\Omega_1)$ of 
 two flat contravariant metrics $\O_2$ and $\O_1$ on $M$ satisfying 
 \begin{enumerate}
     \item $\O_2+\lambda \O_1$ defines a flat metric on $T^*M$ for  a generic constant $\lambda$,
     \item the Christoffel symbols of $\O_2+\lambda \O_1$ are   $\Gamma_{2k}^{ij}+\lambda \Gamma_{1k}^{ij}$, where   $\Gamma_{2k}^{ij}$ and $ \Gamma_{1k}^{ij}$ are the Christoffel symbols of $\O_2$ and $\O_1$,  respectively. 
 \end{enumerate}  
\ed  
\bd \label{FPM}  A  flat pencil of metrics $(\O_2,\O_1)$ on  $M$ is called quasihomogeneous flat pencil of metrics (QFPM) of  degree $d$ if there exists a function $\tau$ on $M$ such that the 
 vector fields $E$ and $e$ defined by 
\begin{eqnarray} \label{tau flat pencil} E&=& \nabla_2 \tau, ~~E^i
=\O_2^{ij}(u)\partial_{u^j}\tau
\\\nonumber  e&=&\nabla_1 \tau, ~~e^i
= \O_1^{ij}(u)\partial_{u^j}\tau  \end{eqnarray} satisfy 
\be \label{vector fields} [e,E]=e,~~ \Lie_E \O_2 =(d-1) \O_2,~~ \Lie_e \O_2 =
\O_1~~\mathrm{and}~~ \Lie_e\O_1
=0.
\ee 
Such a QFPM is  \textbf{regular} if  the
(1,1)-tensor
\begin{equation}\label{regcond}
  R_i^j = \frac{d-1}{2}\delta_i^j + {\nabla_1}_i
E^j
\end{equation}
is  nondegenerate on $M$.
\ed

A Dubrovin-Frobenius manifold is a manifold  with
a smooth structure of a Frobenius algebra on the tangent space  at any point  with certain compatibility conditions. By a Frobenius algebra, we mean  a commutative associative algebra with identity $e$ and a nondegenerate bilinear form $\Pi$ invariant under the product, i.e., $\Pi(a\cdot b,c)=\Pi(a,b\cdot c)$. Globally, we  require the metric $\Pi$ to be flat and the identity vector field $e$ is constant  with respect to it. Detailed information about Dubrovin-Frobenius manifolds and related topics can be found in $\cite{DuRev}$.

Let $M$ be a Dubrovin-Frobenius manifold. In the flat  coordinates $(t^1,...,t^r)$ for $\Pi$ where $e= \partial_{t^{r}}$, the compatibility conditions imply  that there exists a function $\mathbb{F}(t^1,...,t^r)$ which encodes  the Dubrovin-Frobenius manifold structure, i.e., the flat metric is given by
\be \label{flat metric} \Pi_{ij}(t)=\Pi(\partial_{t^i},\partial_{t^j})=  \partial_{t^{r}}
\partial_{t^i}
\partial_{t^j} \mathbb{F}(t)\ee
and, setting $\Omega_1(t)$ to be the inverse of the matrix $\Pi(t)$, the structure constants of the Frobenius algebra are 
\[ C_{ij}^k(t)=\Omega_1^{kp}(t)  \partial_{t^p}\partial_{t^i}\partial_{t^j} \mathbb{F}(t)\]
 The associativity of Frobenius algebra implies $\mathbb F(t)$ satisfies the WDVV equations
 \begin{equation} \label{frob}
 \partial_{t^i}
\partial_{t^j}
\partial_{t^k} \mathbb{F}(t)~ \Omega_1^{kp} ~\partial_{t^p}
\partial_{t^q}
\partial_{t^n} \mathbb{F}(t) = \partial_{t^n}
\partial_{t^j}
\partial_{t^k} \mathbb{F}(t) ~\Omega_1^{kp}~\partial_{t^p}
\partial_{t^q}
\partial_{t^i} \mathbb{F}(t),~~ \forall i,j,q,n. 
  \end{equation}
 The  definition of Dubrovin-Frobenius manifolds includes the existence  of an Euler vector field $E$ of the form $E=(a_i^j t^i+b^j)\partial_{t^j}$ satisfying
\[ 
E\mathbb F(t)= \left(3-d \right) \mathbb{F}(t)+ \frac{1}{2}A_{ij} t^i t^j+B_i t^i+c.
\]
where  $a_i^j$, $b_j$, $c$, $A_{ij}$ and $B_i$ are constants with $a_r^r\neq 0$. Dubrovin-Frobenius manifold is called algebraic (resp. logarithmic) if the potential is algebraic (resp. logarithmic). Moreover, there is a quasihomogenius flat pencil of metric (not necessary regular) of degree  $d$ associated to the Frobenius structure on $M$ consists of the contravariant metrics, the intersection form ${\Omega}_2^{ij}$ and the flat metrics $\O_1^{ij}$, where 
\begin{eqnarray}
   \O^{ij}_2&:=& E^p(t)\eta^{ik}\eta^{jm} \partial_{t^m}\partial_{t^k}\partial_{t^p}\mathbb F(t). 
\end{eqnarray}
On the other hand,  any regular quasihomogenous flat pencil of metrics  on $M$ defines a unique Dubrovin-Frobenius manifold structure on $M$ \cite{DFP}.

\section{Drinfeld-Sokolov bihamiltonian structure}

Let $\g$ be a complex simple Lie algebra of rank $r$ with the Lie bracket $[\cdot,\cdot]$. Define the adjoint representation $\ad: \g \to \textrm{End}(\g)$ by $\ad_{g_1}(g_2):=[g_1,g_2]$. For $g\in \g$, let  $\mathcal O_g$ denote  the orbit of $g$ under the adjoint group action and $\g^g$ denote the centralizer of $g$ in $\g$, i.e. $\g^g:=\ker \ad_g$. An element  $g$ is called nilpotent if $\ad_g$ is  nilpotent in $\textrm{End}(\g)$. A nilpotent element $g$ is called regular (resp. subregular) if $\dim \g^g=r$ (resp. $\dim \g^g=r+2$).  The set of all regular (resp. subregular) nilpotent elements in $\g$ form one complete adjoint orbit. 

We fix a nilpotent element   $L_1$ in $\g$. The following construction gives  same results if we choose another representative  from    $\mathcal{O}_{L_1}$ (for more  details see \cite{mypaper4} and the references within).  Using Jacobson-Morozov theorem, we fix a nilpotent element $f$ and a semisimple element $h$ such that $A:=\{L_1,h,f\}\subseteq \g$ is  a $sl_2$-triple with relations \begin{equation} \label{sl2:relation1} [h,L_1]= L_1,\quad [h,f]=-
f,\quad [L_1,f]= 2 h.
\end{equation}
We normalize the Killing form on $\g$ to get an invariant bilinear form $\bil . . $ such that $\bil {L_1} f=1$. The eigenvalues of $\ad_h$ are integers and half integers and they lead to Dynkin grading
\begin{equation}\label{grad1}
\g=\bigoplus_{i\in \frac{\mathbb{Z}}{2}} \g_i;~~~~\g_i:=\{g\in\g: \ad_h g= i g\}.
\end{equation}
Let    $\eta_r$ denote the maximal eigenvalue of $\ad_h$.  We  fix  an element $K_1\in \g_{-\eta_r}$.  Let $n:=\dim \g^f$ and $\eta_{1},\ldots,\eta_{n}$ are the eigenvalues of $\ad_h$ on $\g^{L_1}$. 

 We consider the loop algebra $\lop \g$, i.e., the space of smooth functions from the unit circle to $\g$. We extend $\bil . .$ on $\g$ to $\lop \g$ by setting
\begin{equation} (g_1|g_2)=\int_{S^1}\bil {g_1(x)}{g_2(x)} dx;~~~ g_1,g_2 \in \lop \g.
\end{equation}
Then we define the gradient $\delta \f (g)$ for a functional $\f$ on $\lop\g$
 to be the unique element in
$\lop\g$ satisfying
\begin{equation}
\frac{d}{d\theta}\f(g+\theta
{{\mathrm{w}}})\mid_{\theta=0}=\big(\delta \f(g)|{\mathrm{w}}\big)
~~~\textrm{for all } {\mathrm{w}}\in \lop\g.
\end{equation}
This allows us to introduce the following standard compatible local Poisson  brackets  $\{.,.\}_2$ and $\{.,.\}_1$ on $\lop \g$: For any functionals $\I$ and $\f$ on $\lop \g$
\begin{eqnarray}\label{Poissbrac}
\{\f,\I\}_2(g(x)) &:=&\Big( {[K_1,\delta \I(g(x))]}\Big |{\delta
\f(g(x))}\Big), \\\nonumber
\{\f,\I\}_1(g(x)) &:=& \Big( {\partial_x\delta \I(g(x))+[g(x),\delta \I(g(x))]}\Big |{ \delta
\f(g(x))}\Big).
 \end{eqnarray}
We denote their Poisson structures or tensors by  $\mathbb B_2$ and  $\mathbb B_1$, respectively. The Poisson brackets $\mathbb B_2$ and $\mathbb B_1$ form an exact Poisson pencil, i.e., there exists a vector field $V$ called Liouville vector field such that $\Lie_V \mathbb B_2=\mathbb B_1$ and $\Lie_V \mathbb B_1=0$.

We fix  Slodowy slice $Q:=L_1+\g^f$ as a transverse subspace to $\mathcal O_{L_1}$ at $L_1$. We introduce the  coordinates $(z^1,\ldots, z^n)$ for $Q$ such that
\be\label{coord}
Q=L_1+ \sum z^i X_i,~, X_i\in \g^f, \ad_h X_i=- \eta_i X_i, i=1,\ldots, n.
\ee
where $X_r=K_1$. We assign $\deg z^i=\eta_i+1$ and we consider the affine loop space $\Q:=L_1+\lop{\g^f}$

\bt
\cite{mypaper4}\label{psred}
The space $\Q$ inherits a  compatible local Poisson structures   $\mathbb B_2^\Q$ and  $\mathbb B_1^\Q$ from  $\mathbb B_2$ and  $\mathbb B_1$,  respectively. They  can be obtained equivalently by using the bihamiltonian reduction with Poisson tensor procedure, Dirac reduction and the generalized  Drinfeld-Sokolov reduction.
\et

 Let us summarize the Poisson tensor procedure to compute the reduced bihamiltonian structure. Let $z\in \Q \cong \lop{\g^f}$ and  $w\in T_z^*\Q$. We identify $T_z^*\Q$ with $\lop {\g^{L_1}}$ using the bilinear from $\bil . .$.  Then there is a unique lift  $v\in T_z^*\lop\g$  of $w$ satisfying
$
\mathbb B_2(v) \in T_z\Q$. The reduced Poisson structure $\mathbb B_2^\Q$ in this case  is given by
\begin{equation}
\mathbb B_2^Q(w)=\mathbb B_2(v).
\end{equation}
We extend the basis $X_i$ on $\g^f$ to a basis on all of $\g$. The local Poisson bracket $\mathbb B_2^Q$ is known as the classical $W$-algebra associated to $L_1$ (see \cite{fehercomp} for the definition) and form with $\mathbb B_1^Q$ an exact Poisson pencil with respect to Liouville  vector field $\partial_{z^r}$.  The Poisson brackets are finite summation of terms in the form  \cite{DZ}
\begin{eqnarray}\label{loc poiss pre}
 \{z^i(x),z^j(y)\}_m & =& \sum_{k\geq -1}  \{z^i(x),z^j(y)\}^{[k]}_m ; ~m=1,2;~~i,j=1,\ldots,n,\\\nonumber
  \{z^i(x),z^j(y)\}^{[-1]}_m &=& F^{ij}_m(z(x))\delta(x-y) \\\nonumber
  \{z^i(x),z^j(y)\}^{[0]}_m &=& \Omega^{ij}_m(z(x)) \delta' (x-y)+ \Gamma_{k,m}^{ij}(z(x)) z_x^k \delta (x-y),  \\\nonumber
  \{z^i(x),z^j(x)\}^{[k]}_m & = & S^{ij}_{m;k}(z(x)) \delta^{k+1}(x-y)+\ldots,~  ~k>0.
\end{eqnarray}
where  $\delta(x-y)$ is the Dirac  delta function defined by
$\int_{S^1} f(y) \delta(x-y) dy=f(x)$.
We recall that such local bihamiltonian structure  admits a dispersionless limit if  the matrices $F^{ij}_m=0,~m=1,2$. Moreover, if it admits a dimensionless limit then the local Poisson brackets  $\{z^i(x),z^j(y)\}^{[0]}_m$ define a local bihamiltonian structure of hydrodynamic type. Moreover, if the matrices $\Omega^{ij}_m$ are nondegenerate  then the compatibility  of Poisson structures  implies that they define a flat pencil of metrics  \cite{DN}. This gives a connection between local bihamiltonian structures and flat pencil of metrics used to construct examples of Dubrovin-Frobenius manifolds. 

\subsection{Common equilibrium points}

We identify Slodowy slice $Q$ with the subspace of constant loops of $\Q$. Then   $F^{ij}_2(z)$ and $F^{ij}_1(z)$ define compatible Poisson structures $B_2^Q$ and $B_1^Q$, respectively, on $Q$ (a bihamiltonian structure). Moreover, $B_2^Q$ is the transverse Poisson structure of the Lie-Poisson structure on $\g$ \cite{mypaper4}.

Using  Chevalley's theorem, we fix  a complete  system of homogeneous generators $P_{1},\ldots,P_{r}$ of the  ring of invariant polynomials under the adjoint group action.
Let $\overline P_i^0$ denotes the restriction of the invariant polynomial $P_i$ to $Q$.
Then
 $\overline P_1^0(q+\lambda K_1),\ldots, \overline P_r^0(q+\lambda K_1)$  form a complete set of independent  Casimirs of the Poisson pencil $B_\lambda^Q:=B_2^Q+ \lambda B_1^Q$ for all  $\lambda\in \mathbb C$. Following  the argument shift method \cite{bolv1}, we consider  the family of functions
\be \label{define F}
\mathbf{F}:=\cup_{\lambda\in \overline \mathbb{C}}\{P'_\lambda: P'_\lambda ~\mathrm{is~ a~ Casimir~ of} ~B_\lambda^Q\}.
\ee
This family commutes pairwise with respect to both Poisson brackets. The main purpose for applying argument shift method is to show that $\mathbf F$ contains enough number  of functionally independent  functions in order to get a completely integrable system for $B_2^Q$. We explored this problem in  \cite{mypaper5}.


We are interested on the space of common equilibrium points  $N$ of the family $\mathbf F$ 
 \be  N:=\{q\in Q: B_\lambda^Q(dP')(q)=0,~ \forall P'\in \mathbf F, \lambda \in \overline{\mathbb C}\} \ee 
which was introduced and studied in \cite{bolv1} for arbitrary finite dimensional bihamiltonian structure. 
 
\section{Nilpotent  elements of regular semisimple type} \label{regular type}

From this section and the following sections, we suppose   $L_1$ is a nilpotent element of regular semisimple type. The classification of such elements  were obtained in \cite{Elash}. By definition, we can and we will fix $K_1$ such that $h':=L_1+K_1$ is regular semisimple element. Thus, $\h':=\ker \ad_{h'}$ is a Cartan subalgebra known as opposite Cartan subalgebra.   The  adjoint group element $w:=\exp {2\pi \mathbf i\over \eta_r+1} \ad_h$  acts on $\h'$ as a representative of regular  conjugacy class of order $\eta_r+1$ or $2 \eta_r+2$ in the underline  Weyl group $\W(\g)$. We can and we will arrange the basis  $L_j's$   such that for $i=1,\ldots, r$, there exists $K_i \in \oplus_{i\leq 0}\g_i$ such that  $Y_i=L_i+K_i\in \h'$ and $Y_1,\ldots,Y_r$ from a basis  for $\h'$. Moreover, we can and we will arrange the basis $X_i$ of $\g^f$ such that $\bil {X_i}{L_j}=\delta_{ij}$ and we  assume $0\leq \eta_1\leq \ldots\leq \eta_r$.  



For comparison and fixing notations we summarize briefly the results of \cite{mypaper6}. We say the  element $L_1$ is  distinguished if the orbit  $\mathcal O_{L_1}$  has no representative in a proper Levi subalgebra of $\g$. If $L_1$ is a distinguished nilpotent element of regular semisimple type then  $w$ will be a representative of  regular cuspidal conjugacy class \cite{DelFeher}. This relation gives  a correspondence between  distinguished nilpotent orbits in simple Lie algebras   and regular cuspidal conjugacy classes in  Weyl groups. Moreover, the eigenvalues of $w$ has the form $\exp \frac{2\pi \mathbf{i}\eta_i}{\eta_r+1}$, $i=1,\ldots,r$.

In the remaining of this section, we assume $L_1$ is a distinguished nilpotent element of regular semisimple type. Then, it turns out that there is a quasihomogeneous  change of coordinates on $Q$ in the form 
 \be
 t^i=
 \left\{
  \begin{array}{ll}
z^1, & \hbox{i=1} , \\
  z^i+\mathrm{non~ linear~ terms},&\hbox{i=2,\ldots,r,} \\
z^i,&\hbox{i=r+1,\ldots,n.}
  \end{array}
\right.
\ee
such that 
 \begin{enumerate}
      \item $\deg t^i=\deg z^i=\eta_i+1$
     \item $t^1,\ldots,t^r$ form a complete set of  Casimirs of $B_1^Q$ and they are in involution with respect to $B_2^Q$.
     \item For $1\leq i\leq r$, $t^i$  is a linear  combination  of invariant polynomials  $\overline P_j^0$ and their derivative with respect to  ${z^r}$.
 \end{enumerate}

The special coordinates $(t^1,\ldots,t^r)$ are the main ingredient to prove the following theorems.
\bt\label{rank of B1} \cite{mypaper6}
The family $\mathbf F$ is complete (contains completely integrable system for $B_2^Q$) for every distinguished nilpotent element of semisimple type. In particular, $\rk\,B_1^Q=n-r$ and 
 \be\label{def of N}
 N=\{q\in Q: \ker B_1^Q(q)=\ker B_2^Q(q)\}.
 \ee
\et

\bt \label{main thm} \cite{mypaper6}
Let $\g$ be a complex simple Lie algebra of rank $r$. Fix a regular cuspidal conjugacy class  $[w]$  in the Weyl group $\W (\g)$ of $\g$. Assume the order of representatives in $[w]$ is  $\eta_r+1$ and eigenvalues are $\epsilon^{\eta_i}$, $i=1,\ldots,r$, where $\epsilon$ is a primitive $(\eta_r+1)$th root of unity. Let $\mathcal{O}_{L_1}$ be the distinguished  nilpotent orbit of semisimple type  associated to $[w]$  under the notion of opposite Cartan subalgebra. Consider  the finite bihamiltonian structure formed by the leading term of Drinfled-Sokolov bihamiltonain structure associated to a representative $L_1$ of $\mathcal{O}_{L_1}$. Then its space of common equilibrium points   acquires an algebraic Dubrovin-Frobenius manifold structure with charge $\frac{\eta_r-1}{\eta_r+1}$ and degrees $\frac{\eta_i+1}{\eta_r+1}$. This structure depends only on the conjugacy class. 
\et

Let $s$ be the number of zero notes in the weighted Dynkin Diagram associated to $L_1$. Then, after  normalization, the space $N$ is defined by 
 \begin{eqnarray}
 N &=& \{t: \partial_{t^\beta} \overline P_j^0(t)=0; j=r-s+1,\ldots,r,~~\beta=r+1,\ldots,n\}.\label{defEqN1}
 \end{eqnarray}
Moreover, $(t^1,\ldots,t^r)$ provide local coordinates around generic points of  $N$. 
Let $\N=\lop N$ be the loop space of $N$.  Then Dirac reduction of the Poisson pencil  $\mathbb B_\lambda^\Q:=\mathbb B_2^\Q+\lambda \mathbb B_1^\Q$ to $\N$ is well defined and leads to compatible local Poisson brackets $\{.,.\}_\alpha^\N$, $\alpha=1,2$ which admit  a dispersionless limit and form an exact Poisson pencil with respect to the vector field $\partial_{t^r}$. Moreover, $\{.,.\}_2^\N$ is a  classical $W$-algebra.

The local Poisson brackets on $\N$  are simply given  by taking the upper left minors of  Poisson brackets \eqref{loc poiss pre}  under the coordinates  $(t^1,\ldots,t^n)$, i.e.,
\begin{eqnarray}\label{loc poiss}
  \{t^a(x),t^b(y)\}^{[-1]}_\alpha &=& 0,~ a,b=1,\ldots r, ~\alpha=1,2.\\\nonumber
  \{t^a(x),t^b(y)\}^{[0]}_\alpha &=& \O^{ab}_\alpha(t(x)) \delta' (x-y)+ \Gamma_{\alpha k}^{ab}(t(x)) t_x^k \delta (x-y)\\\nonumber   \{t^a(x),t^b(x)\}^{[k]}_\alpha & = & S^{ab}_{\alpha;k}(t(x)) \delta^{k+1}(x-y)+\ldots,~  ~k>0.
\end{eqnarray}
where $t^k, k>r$ are  solutions  of the polynomial equations \eqref{defEqN1} defining $N$. Then the non-degeneracy   of the matrices $\O^{ab}_\alpha(t)$ is deduced from the non-degeneracy  of the restriction of the bilinear from $\bil . .$ to $\h'$. In addition,  they  define a regular quasihomogeneous  flat pencil of metrics on $N$ of degree $d=\frac{\eta_n-1}{\eta_n+1}$ which  leads to an algebraic Dubrovin-Frobenius structure. 

\subsection{Central invariants}

We recall from \cite{DLZ} the following. Suppose that the bihamiltonian structure has the form \eqref{loc poiss} and the flat pencil of metrics defined by $\O_{2}$ and $\O_{1}$  leads to a Dubrovin-Frobenius manifold structure.  Then the Dubrovin-Frobenius manifold structure is called  semisimple if the roots $u^1,\ldots,u^r$ of the characteristic polynomial
\be \label{semsimpfp}
\Psi(\lambda;t):=\det (\Omega^{uv}_2(t)-\lambda \Omega^{uv}_1(t))
\ee 
are pairwise distinct at generic points on $M$. In this case  $(u^1,\ldots,u^r)$ define local coordinates.  Moreover, we can define the central invariants $c^1(u^1),\ldots,c^r(u^r)$ of the local bihamiltonian structure under generalized Miura type transformations  \cite{DLZ}. When $ \{t^a(x),t^b(x)\}^{[1]}_\alpha=0$, the central invariants  can be calculated by the formulas
\be \label{cent form}
c_i(u^i):=\frac{1}{3}[\frac{d\Psi}{d\lambda}(u^i;t)]^2\frac{\frac{\partial \Psi(\lambda;t)}{\partial t^k}\frac{\partial \Psi(\lambda;t)}{\partial t^l}(S_{2;2}^{kl}(t)-\lambda S_{1;2}^{kl}(t))}{[\frac{\partial \Psi(\lambda;t)}{\partial t^k}\frac{\partial \Psi(\lambda;t)}{\partial t^l}\Omega_1^{kl}(t)]^2}{\Big |}_{\lambda=u^i}.
\ee 
We recall that changing  representatives of the bihamiltonian structure  by   linear combinations leads to a fractional linear transformation for the central invariants. In particular, under rescaling by a constant $\kappa$, i.e., 
 \be
\{.,.\}_\alpha\mapsto \kappa  
\{.,.\}_\alpha, ~~\alpha=1,2.
 \ee
 we get $c^i\mapsto \kappa^{-1} c^i$ and an equivalent Frobenius manifold structure \cite{DuRev}. Note that if the local bihamiltonian structure forms an exact Poisson pencil then the central invariants are constant \cite{FalLor}.
 
 Suppose the central invariant are all equal and constant. Then  the bihamiltonian structure is of of topological type, i.e., it is the bihamiltonian structure associated with the principle hierarchy of the Dubrovin-Frobenius manifold structure \cite{DZ}.  Note that under  suitable rescaling  the central invariants all equaling $\frac{1}{24}$ which is the number obtained in \cite{DZ}.

Assume the nilpotent element $L_1$ is regular.  Then we get the standard Drinfeld-Sokolov bihamiltonian structure \cite{DS}. Here, $Q=N$ and  theorem \ref{main thm} leads to polynomial Dubrovin-Frobenius manifold structure. In this case the central invariants are calculated in \cite{DLZ}. It turns out that the bihamiltonian structure is of topological type only when the Lie algebra $\g$ is simply-laced.  The central invariants of the bihamiltonian structures associated to algebraic non-polynomial Dubrovin-Frobenius manifolds obtained by theorem \ref{main thm} have not been  computed yet. 

 Information about the notion of central invariants  and its applications can be found in \cite{DLiZ}, \cite{LZ1} and \cite{Lo}.

\section{Lie algebra $A_3$ and $A_4$}

In the reminder of this article, we follow the steps of section \ref{regular type}  for  the non-distinguished nilpotent elements of regular semisimple type in Lie algebra $sl_3$ or $sl_4$ which are simple Lie algebras of type $A_3$ and $A_4$, respectively. Precisely, in the classification of nilpotent orbits in $sl_n$ by the integer partition of $n$, we assume $L_1$ is a representative of the nilpotent orbit    corresponding  to  the partition [2,2], [3,1] in $sl_4$, or  [2,1] in $sl_3$ \cite{Elash}. The partition of subregular nilpotent orbit  has the form  $[n-1,1]$. Let $\varepsilon_{i,j}$ denote the basic square matrix  with 1 on the $(i,j)$-entry.

\subsection{Partition $[2,1]$ in $sl_3$}

We take  $L_1:=\varepsilon _{1,3}$ as a representative of the nilpotent orbit   $[2,1]$ in $sl_3$.  We fix the $sl_2$-triple 
\be 
L_1=\varepsilon _{1,3}, ~ h=\frac{1}{2} \left(\varepsilon
   _{1,1}-\varepsilon
   _{3,3}\right), ~ f=\varepsilon _{3,1}.
\ee
The numbers $\eta_i's$ in the standard order given in section \ref{regular type} are $1,0,\frac{1}{2},\frac{1}{2}$.  Elements of Slodowy slice will have the form 
\be 
q=L_1+\sum_{i=1}^4 z_i X_i=\left(
\begin{array}{ccc}
 z_1 & 0 & 1 \\
 z_3 & -2 z_1 & 0 \\
 z_2 & z_4 & z_1
\end{array}
\right).
\ee 
 For the restriction of invariant polynomials  to $Q$, we fix 
\be 
\overline{P}_1^0=z_2+3 z^2_1,~ \overline{P}_2^0=z_2 z_1- z_1^3+\frac{1}{2}z_3 z_4.
\ee 

Let $\I$ be a functional  on $Q$  and $\delta \I=(\frac{\delta \I}{\delta z_1},\frac{\delta \I}{\delta z_2},\frac{\delta \I}{\delta z_3},\frac{\delta \I}{\delta z_4})$. Using Poisson tensor procedure, the  lift $\overline{\delta \I}$  of  $\delta \I$ in $\lop \g$ takes the form 
\be 
\overline{\delta \I}=\left(
\begin{array}{ccc}
\frac{1}{6}\frac{\delta \I}{\delta z_1}+v_1 &  \frac{\delta \I}{\delta z_3} &  \frac{\delta \I}{\delta z_2} \\
 v_3 & -\frac{1}{3} \frac{\delta \I}{\delta z_1} &  \frac{\delta \I}{\delta z_4}\\
 v_2 & v_4 & \frac{1}{6} \frac{\delta \I}{\delta z_1}-v_1
\end{array}
\right)
\ee 
where the values of $v_1,\ldots,v_4$ are uniquely determined by the requirement  
\be \partial_x \overline{\delta \I}+[q(x),\overline{\delta \I}]\in T_q\Q\cong \lop {\g^f}. \ee
Then the local Poisson bracket $\mathbb B_2^Q$ for two functional $\f$ and $\I$ on $\Q$ is given by 
\be  \{\f[q(x)],\I[q(y)]\}_2^{L_1}=\Big<\overline{\delta \f}| \frac{d}{dx}(\overline{\delta \I})+[q(x),\overline{\delta \I}]\Big>.  \ee In coordinates, $\mathbb B_2^Q$ (the classical $W$-algebra) has the following nonzero brackets. 
\begin{eqnarray}
\{z_1(x) ,z_1(y)\}&=&   \frac{1}{6} \delta',~~
   \{z_1(x) ,z_3(y)\}=-\frac{1}{2} z_3 \delta,~~\{z_1(x) ,z_4(y)\}=\frac{1}{2} z_4 \delta  \\\nonumber
   \{z_2(x) ,z_2(y)\}&=& 
   z_2'\delta+2 z_2
   \delta '-\frac{1}{2}
   \delta ^{(3)},~    \{z_2(x) ,z_3(y)\}=\frac{1}{2} \left(
   \left(z_3'+6
   z_1 z_3\right)\delta +3 z_3
   \delta '\right)
  \\\nonumber
   \{z_2(x) ,z_4(y)\}&=&   \frac{1}{2} \left( 
   \left(z_4'-6
   z_1(x) z_4\right)\delta+3 z_4
   \delta \right),~\{z_2(x) ,z_4(y)\} =\frac{1}{2} z_4 \delta   \\\nonumber \{z_3(x) ,z_4(y)\}&=&- \left(3
   z_1'-9
   z_1^2+z_2\right)\delta -6
   z_1 \delta '+\delta
   ''.
\end{eqnarray} 
Here and in what follows, for local Poisson brackets all function on the right hand side are functions of $x$ and we write $\delta$ instead of $\delta(x-y)$.

 There are two bihamiltonian structures in the literature associated to this classical $W$-algebra: One is constructed using the above mentioned theory of opposite Cartan subalgebra and another one is associated to  the fractional KdV hierarchy \cite{gDSh2}. We treat them in details below.
\subsubsection{Using opposite Cartan subalgebra}
The opposite Cartan subalgebra has basis $X_1$ and $L_1+f$.
 Then Drinfeld-Sokolov bihamiltonian structure is formed by $\mathbb B_2^Q$ and   $\mathbb B_1^\Q=\Lie_{\partial_{z_2}}\mathbb B_2^\Q$. In the notations of \eqref{loc poiss pre}, we have  \be 
F_2^{ij}(z)=\left(
\begin{array}{cccc}
 0 & 0 & -\frac{1}{2} z_3&
   \frac{1}{2}z_4 \\
 0 & 0 & 3 z_1 z_3 & -3 z_1 z_4 \\
 \frac{1}{2}z_3 & -3 z_1 z_3 & 0 &
   9 z_1^2-z_2 \\
 -\frac{1}{2}z_4 & 3 z_1 z_4 &
   z_2-9 z_1^2 & 0
\end{array}
\right)
,~~
\O_2^{ij}(z
)=\left(
\begin{array}{cccc}
 \frac{1}{6} & 0 & 0 & 0 \\
 0 & 2 z_2 & \frac{3}{2} z_3 &
   \frac{3}{2} z_4 \\
 0 & \frac{3}{2}  z_3& 0 & -6 z_1
   \\
 0 & \frac{3}{2} z_4 & -6 z_1 & 0
\end{array}
\right)
\ee 
The space of common equilibrium points $N$ is defined by
\be
N:=\{q\in Q: \partial_{z_3}\overline P_2^0=0=\partial_{z_4}\overline P_2^0\}=\{q\in Q: z_3=0=z_4\}
\ee

The special coordinates $(t_1,\ldots,t_4)$ are defined by 
\be 
t_1=\overline P_1^0,~ t_2=\partial_{z_2} \overline P_2^0=z_1,~ t_3=z_3,~t_4=z_4.
\ee 
Then $(t_1,t_2)$ form a local coordinates on $N$ and the reduced bihamiltonian  structures  on  the loops space $\N:=\lop N$ using Dirac reduction has the  leading terms. 
\be 
\O_2^{ab}(t)=\left(
\begin{array}{cc}
 2 t_1 & t_2 \\
 t_2 & \frac{1}{6}
\end{array}
\right), S_{2;1}^{ab}=0, S_{2;2}^{ab}(t)=\left(
\begin{array}{cc}
 -\frac{1}{2} & 0 \\
 0 & 0
\end{array}
\right)
\ee 
In contrast with the case of distinguished nilpotent elements, $\O_1^{ab}=\Lie_{\partial_{t_1}} \O_2^{ab}$ degenerate and does not define a metric on $N$. However,  $\Lie_{\partial_{t_2}}^2 \O_2^{ab}=0$ and $\Lie_{\partial_{t_2}} \O_2^{ab}$ is nondegenerate. It turns out that  $\Lie_{\partial_{t_2}}^2 \mathbb B_2^\N=0$. Thus $\mathbb B_2^\N$ and $\overline \mathbb B_1^\N:=\Lie_{\partial_{t_2}} \mathbb B_2^\N=0$ form a new bihamiltonian structure on $\N$ \cite{Ser}. Moreover, the matrices  $\O_2^{ab}$ and  $\Lie_{\partial_{t_2}} \O_2^{ab}=0$ from a regular quasihomogeneous flat pencil of metric on $N$. The corresponding  Dubrovin-Frobenius manifold structure on $N$ is equivalent to the Dubrovin-Frobenius manifold structure constructed on the  Hurwitz space $M_{0;1,0}$. It  has the potential
\be 
\mathbb F(t_1,t_2)=\frac{1}{2} t_2^2 t_1+\frac{1}{24}
   t_1^2 \log t_1.
\ee 
The Euler vector field is $E=2 t_1\partial_{t_1}+t_2\partial_{t_2}$ with charge $d=-1$. Note that $E\mathbb F=(3-d)\mathbb F+\frac{1}{12} t_1^2$. Using the formula \eqref{cent form}, we compute  the  central invariants of the corresponding bihamiltonian structure and they equal $-\frac{1}{24}$. Thus the bihamiltonian structure is of topological type. 

\subsubsection{For Fractional KdV Hierarchy}\label{frac}
We mention that one can get more bihamiltonian structure on $\Q$ for non-distinguished nilpotent elements for certain choices of $K_1$. 
At the definition of $B_1$ in \eqref{Poissbrac}, we fix   $K_1=\varepsilon_{2,1}+\varepsilon_{2,3}$ instead of $\varepsilon_{3,1}$. Then we get a local Poisson bracket  $\mathbb B_1^Q$  compatible with $\mathbb B_2^Q$ and they form an exact Poisson pencil using (with minor modification) the reductions mentioned in theorem \ref{psred} \cite{mypaper4}. To find $ \mathbb B_1^Q$, we  simply   introduce the change of coordinates 
\[ t_1= 2 z_1^2+\frac{2
   z_2}{3},t_2=z_3+z_4, t_3=z_3-z_4, t_4=z_1.\]
then  \[\mathbb B_1^Q=\Lie_{\partial
_{t_2}} \mathbb B_2^Q.\] 
Here, the restrictions of the invariant polynomials will be \[ \overline P_1^0(t)=t_1, ~\overline P_2^0(t)=-8 t_4^3+3 t_1
   t_4+\frac{t_2^2}{4}-\frac{t_3^2}{4}\]
and 
\[ F^{ij}_2(t)=\left(
\begin{array}{cccc}
 0 & 0 & 0 & 0 \\
 0 & 0 & 3 t_1-24 t_4^2 & \frac{t_3}{2}
   \\
 0 & 24 t_4^2-3 t_1 & 0 & \frac{t_2}{2}
   \\
 0 & -\frac{t_3}{2} & -\frac{t_2}{2} &
   0
\end{array}
\right),~\Omega_2^{ij}(t)= \left(
\begin{array}{cccc}
 \frac{4 t_1}{3} & t_2 & t_3 & \frac{2
   t_4}{3} \\
 t_2 & -12 t_4 & 0 & 0 \\
 t_3 & 0 & 12 t_4 & 0 \\
 \frac{2 t_4}{3} & 0 & 0 & \frac{1}{6}
\end{array}
\right)\]
space of common equilibrium points $N$ is defined by
\be
N:=\{q\in Q: \partial_{t_3}\overline P_2^0=0=\partial_{t_4}\overline P_2^0\}=\{q\in Q: t_3=0,~t_4=\sqrt{\frac{t^1}{8} }\}.
\ee
Here $(t_1,t_2)$ form local coordinates on $N$ and the reduced bihamiltonian  structures on $\lop N$ has the  leading terms. 
\be 
\O_2^{ab}(t)=\left(
\begin{array}{cc}
 \frac{4 t_1}{3} & t_2 \\
 t_2 & 3 \sqrt{2} \sqrt{t_1}
\end{array}
\right),~ S_{2;1}^{ab}=0, S_{2;2}^{ab}(t)=\left(
\begin{array}{cc}
 -\frac{2}{9} & 0 \\
 0 & 0
\end{array}
\right)
\ee 
Note that $\O_1^{ab}=\Lie_{\partial_{t_2}} \O_2^{ab}$ from with $\O_2^{ab}$ a regular quasihomogenius flat pencil of metrics. 
 The corresponding  Dubrovin-Frobenius manifold structure on $N$ has the potential 
\be 
\mathbb F(t_1,t_2)=\frac{6}{5} \sqrt{2}
   t_1^{5/2}+\frac{1}{2} t_2^2 t_1
\ee 
The Euler vector field is $E=\frac{4}{3} t_1\partial_{t_1}+t_2\partial_{t_2}$ with charge $d=-\frac{1}{3}$. We calculate the  central invariants and they are all equal $-\frac{1}{54}$ and the bihamiltonian structure on $\N$ is of topological type.

\subsection{Partition $[3,1]$ in $sl_4$}

We fix $L_1=\varepsilon _{1,2}+\varepsilon _{2,3}$ as a representative of the nilpotent orbit of type $[3,1]$ in the Lie algebra $sl_4$.  We introduce  the $sl_2$-triple
\[ L_1=\varepsilon _{1,2}+\varepsilon _{2,3},~ h= \varepsilon _{1,1}- \varepsilon
   _{3,3},~ f=2 \varepsilon _{2,1}+2 \varepsilon
   _{3,2}
   \]
Then $\dim \g^f=5$ and the numbers $\eta_1,\ldots,\eta_5$ are $1,0,2,1,1$, respectively.
We write elements of Slodowy slice in the form 
\[Q=L_1+\sum_{i=1}^5 z_i X_i=\left(
\begin{array}{cccc}
 \frac{1}{\sqrt{3}}z_2 & 1 & 0 & 0 \\
 2 z_1 & \frac{1}{\sqrt{3}}z_2 & 1 & 0 \\
 4 z_3 & 2 z_1 & \frac{1}{\sqrt{3}}z_2 & 2 z_4-2 z_5 \\
 2 z_4+2 z_5 & 0 & 0 & -\sqrt{3} z_2
\end{array}
\right)\]
where $X_i\in \g^f$ and $\ad_h X_i=-2\eta_i$. 
The fix invariant polynomials having the  restrictions   
\begin{eqnarray}
\overline P_1^0&=& z_1+\frac{1}{2} z_2^2,~~\overline P_2^0=z_3-\frac{2}{3 \sqrt{3}} z_2^3+\frac{2}{\sqrt{3}} z_1 z_2\\\nonumber
\overline P_3^0&=& \frac{7}{12} z_2^4+z_1 z_2^2+\sqrt{3} z_3 z_2+2 z_1^2+z_4^2-z_5^2
\end{eqnarray}
The Poisson structure  $\mathbb B_2^\Q$  has the following nonzero brackets. 
\begin{eqnarray}
\{z_1(x) ,z_1(y)\}&=&   z_1' \delta  +2
   z_1  \delta ' -\frac{1}{2}
   \delta ^{(3)},~~
   \{z_1(x) ,z_3(y)\}= 2  z_3' \delta  +3
   z_3  \delta ' \\\nonumber
    \{z_1 (x),z_4(y)\}&=&
   \left(z_4' -\frac{4}{\sqrt{3}}
   z_2  z_5 \right)\delta  +2
   z_4  \delta ',~~
    \{z_1(x) ,z_5(y)\}=
   \left(z_5' -\frac{4}{\sqrt{3}}
   z_2  z_4 \right)\delta  +2
   z_5  \delta '  \\\nonumber
   \{z_2(x) ,z_2(y)\}&=& \delta ' ,~~
     \{z_2(x) ,z_4(y)\}=\frac{4}{\sqrt{3}} z_5  \delta,~~
      \{z_2(x) ,z_5(y)\}=\frac{4}{\sqrt{3}} z_4  \delta   \\\nonumber
              \{z_3(x) ,z_3(y)\}&=&
            \frac{1}{12}\left(-9
   z_1'' +32
   z_1^2 {}+48 z_4^2 {}-48
   z_5^2 {}\right) \delta
   ' + \frac{1}{24}\Big(-30
   z_1'  \delta
   '' -20 z_1
   \delta ^{(3)} +\delta
   ^{(5)} \Big)
     \\\nonumber&+ &            \frac{1}{6} \left(16 z_1
   z_1' +24 z_4
   z_4' - 24
   z_5
   z_5' +z_1
   ^{(3)}{} \right)\delta
     \\\nonumber
   \{z_3 (x),z_4(y)\}&=&-\frac{1}{6}  \Big(4
   \sqrt{3} z_2
   z_4' +12 \sqrt{3}
   z_4
   z_2' -32 z_5  z_2^2 {} -z_5'' +16
   z_1  z_5 \Big)\delta+ \frac{5}{6}
   \left(z_5' -4
   \sqrt{3} z_2  z_4 \right) \delta
   ' +\frac{5}{3} z_5
   \delta ''     \\\nonumber
       \{z_3(x) ,z_5(y)\}&=&  \frac{5}{6}
   \left(z_4' -4
   \sqrt{3} z_2  z_5 \right) \delta
   ' +\frac{5}{3} z_4
   \delta ''-\frac{1}{6}  \Big(4
   \sqrt{3} z_2
  z_5' +12 \sqrt{3}
   z_5
   z_2' -32 z_4  z_2^2 {}-z_4'' +16
   z_1  z_4 \Big)\delta   \\\nonumber
       \{z_4(x) ,z_4(y)\}&=&  -
   \left(z_1' -8 z_2
   z_2' \right)\delta  +\frac{1}{2} \delta
   ^{(3)} - 2
   \left(z_1 -4 z_2^2 {}\right)
   \delta '  \\\nonumber
       \{z_4(x) ,z_5(y)\}&=& \frac{2}{9}  \left(-3
   \sqrt{3} z_2'' -16
   \sqrt{3} z_2^3 {}+12 \sqrt{3}
   z_1  z_2 -9 z_3 \right)\delta  -2
   \sqrt{3} \left(z_2'
   \delta ' +z_2  \delta
   '' \right) \\\nonumber
       \{z_5(x) ,z_5(y)\}&=&
   \left(z_1' -8 z_2
   z_2' \right)\delta  -\frac{1}{2} \delta
   ^{(3)} + 2
   \left(z_1 -4 z_2^2 {}\right)
   \delta '
\end{eqnarray}
Recall that the Poisson structure $\mathbb B_1^\Q$ is  the Lie derivative  of $\mathbb B_2^\Q$ along $\partial_{z_3}$. Using the notation of equations \eqref{loc poiss pre}, the transverse Poisson structure $B_2^Q$ will have the following matrices
\[
F^{ij}_2(z)=\left(
\begin{array}{ccccc}
 0 & 0 & 0 & -\frac{4 z_2
   }{\sqrt{3}}z_5 & -\frac{4}{\sqrt{3}} z_2
   z_4 \\
 0 & 0 & 0 & \frac{4}{\sqrt{3}} z_5 &
   \frac{4}{\sqrt{3}} z_4 \\
 0 & 0 & 0 & \frac{16}{3} z_2^2
   z_5-\frac{8}{3} z_1 z_5 &
   \frac{16}{3} z_2^2 z_4-\frac{8}{3} z_1
   z_4 \\
~* & * & * & 0 &
   -\frac{32}{3 \sqrt{3}} z_2^3+\frac{8}{\sqrt{3}}
   z_1 z_2-2 z_3 \\
~* &* & * &
   \frac{32}{3 \sqrt{3}} z_2^3-\frac{8}{\sqrt{3}}
   z_1 z_2+2 z_3 & 0
\end{array}
\right)
\]
and 
\[ \Omega_2^{ij}(z)=\left(
\begin{array}{ccccc}
 2 z_1 & 0 & 3 z_3 & 2 z_4 & 2 z_5 \\
 0 & 1 & 0 & 0 & 0 \\
 3 z_3 & 0 & \frac{1}{24} \left(64 z_1^2+96 z_4^2-96 z_5^2\right)
   & -\frac{10}{\sqrt{3}} z_2 z_4 & -\frac{10}{\sqrt{3}} z_2 z_5
   \\
 2 z_4 & 0 & -\frac{10}{\sqrt{3}}  z_2 z_4& \frac{1}{4} \left(32
   z_2^2-8 z_1\right) & 0 \\
 2 z_5 & 0 & -\frac{10}{\sqrt{3}}  z_2 z_5& 0 & \frac{1}{4}
   \left(8 z_1-32 z_2^2\right)
\end{array}
\right)\]
and we note that the matrix $S_{2;1}^{ij}$ is not zero.

The space of common equilibrium points $N$ is defined by
\be
N:=\{q\in Q: \partial_{z_4}\overline P_3^0=0=\partial_{z_5}\overline P_3^0\}=\{q\in Q: z_5=0=z_4\}
\ee
Consider the reduced  bihamiltonian structure on loop space $\N$ which admits a dispersionless limit and  has the form \eqref{loc poiss}. The matrix $\Omega^{ab}_2(z)$ is found by taking the upper left minor of $\Omega^{ij}_2(z)$ where  $z_4=z_5=0$. Then  the matrix $\Omega_1^{ab}=\Lie{\partial_{z_3}}\Omega^{ab}_2$  degenerate on $N$. 

We try to find another metric  which forms with $\Omega^{ab}_2$ a flat pencil of metrics. From  conditions \eqref{vector fields}, we seek a vector field $e$ such that $\Lie_{e}\Omega^{ab}_2$ is nondegenerate and $\Lie_{e}^2\Omega^{ab}_2=0$. By  ad-hoc trials, we find the vector field 
\be e= \partial_{z_1}+ 2 i \sqrt{\frac{2}{3}}  z_2 \partial_{z_3}\ee
satisfies the required conditions. It turns out that   $\Lie_{e}^2 \mathbb B_2^\N=0$. Hence, we get a new local bihamiltonian structure on $\N$ formed by $\mathbb B_2^\N$ and $\mathbb \Lie_{e} \mathbb B_2^\N$ \cite{Ser} which  admits a dispersionless limit and form and exact Poisson pencil. Moreover,  the leading term defines a regular quasihomogeneous flat pencil of metrics. In the following  coordinates \be  t_1=z_1+\frac{1}{2}z_2^2,~t_2=i \sqrt{\frac{3}{2}}
   z_2,~ t_3=z_3-\frac{4}{3} i \sqrt{\frac{2}{3}} z_2^3-2 i
   \sqrt{\frac{2}{3}} z_1 z_2\ee 
we get 
\be  \Omega_2^{ab}(t)=\left(
\begin{array}{ccc}
 2 t_1 & t_2 & 3 t_3 \\
 t_2 & -\frac{3}{2} & 2 t_1-2 t_2^2 \\
 3 t_3 & 2 t_1-2 t_2^2 & -8 t_2 t_3
\end{array}
\right), ~~S_{2;2}^{ab}(t)=\left(
\begin{array}{ccc}
 -\frac{1}{2} & 0 & \frac{2 }{3}t_2 \\
 0 & 0 & 0 \\
 \frac{2}{3} t_2 & 0 & -\frac{7 }{6}t_2^2-\frac{5}{6} t_1
\end{array}
\right), ~S_{1;2}^{ab}=0,\ee 
and  $e=\partial_{t_1}$. The resulting Dubrovin-Frobenius manifold structure is equivalent to the Dubrovin-Frobenius structure on Hurwitz space $M_{0;1,0}$ (\cite{DuRev}, example 5.5). It has the logarithmic  potential
\be  \mathbb F=\frac{t_1^3}{12}+\frac{1}{2} t_2 t_3 t_1-\frac{1}{6} t_2^3
   t_3-\frac{1}{8} t_3^2 \log t_3\ee 
   with Euler vector field $E=\frac{1}{2}\sum\Omega^{1i}\partial_{t_i}$ and charge $0$. Note that $E \mathbb F= 3 \mathbb F-\frac{3}{16} t_3^2$. Then using the formula \eqref{cent form}, we find  all central invariants equal $-\frac{1}{96}$.  Thus the corresponding  bihamiltonian structure on $\N$ is of topological type.

\subsection{Nilpotent orbit of type $[2,2]$}
We fix the following $sl_2$-triple where $L_1$ is a representative of the nilpotent orbit of type $[2,2]$ in  $sl_4$. \[ L_1=\varepsilon _{1,2}+\varepsilon _{3,4},~ h=\frac{1}{2}(\varepsilon _{1,1}-\varepsilon
   _{2,2}+\varepsilon
   _{3,3}-\varepsilon _{4,4}),~ f=\varepsilon _{2,1}+\varepsilon
   _{4,3}
   \]
Then the  numbers $\eta_1,\ldots,\eta_7$ are  $1,0,1,1,0,0,1$, respectively. We write elements of Slodowy slice in the form 
\be Q=L_1+\sum_{i=1}^7 z_i X_i= \left(
\begin{array}{cccc}
 \frac{1}{\sqrt{2}}z_2 & 1 &
   \frac{1}{\sqrt{2}}(z_5-z_6) & 0 \\
 z_1+z_3 & \frac{1}{\sqrt{2}}z_2 &
   z_4+z_7 &
   \frac{1}{\sqrt{2}}(z_5-z_6) \\

   \frac{1}{\sqrt{2}}(z_5+z_6) & 0 &
   -\frac{1}{\sqrt{2}}z_2 & 1 \\
 z_7-z_4 &
   \frac{1}{\sqrt{2}}(z_5+z_6) & z_1-z_3 &
   -\frac{1}{\sqrt{2}}z_2
\end{array}
\right)\ee
We fix invariant polynomials having the restrictions  to $Q$ equals 
\begin{eqnarray}
\overline P_1^0 &=& \frac{1}{2}z_2^2+\frac{1}{2}z_5^2-\frac{1}{2}z_6^2+z_1 ,~\overline P_2^0=z_2 z_3+z_4 z_6+z_5 z_7\\\nonumber 
\overline P_3^0&=& \frac{1}{4}z_2^4+\frac{1}{2} z_5^2
   z_2^2-\frac{1}{2} z_6^2 z_2^2-z_1
   z_2^2+\frac{1}{4}z_5^4+\frac{1}{4}z_6^4
   +z_1^2-z_3^2+z_4^2-z_1
   z_5^2-\frac{1}{2} z_5^2 z_6^2+z_1
   z_6^2-z_7^2
\end{eqnarray}
Then the nonzero brackets of  $\mathbb B_2^Q$ are
 \begin{eqnarray}
\{z_1(x),z_1(y)\}&=&-\frac{1}{2}
   \delta ^{(3)}+ z_1'  \delta   +2 z_1  \delta '\\\nonumber
 \{z_1(x),z_3(y)\}&=& \delta
   \left(z_3' -\sqrt{
   2} z_4  z_5 -\sqrt{2} z_6
   z_7 \right)+2 z_3  \delta
   '
   \\\nonumber \{z_1(x),z_4(y)\}&=&
   \left(z_4' -\sqrt{
   2} z_3  z_5 +\sqrt{2} z_2
   z_7 \right)\delta+2 z_4  \delta
   '\\\nonumber \{z_1(x),z_7(y)\}&=&
   \left(z_7' +\sqrt{
   2} z_2  z_4 +\sqrt{2} z_3
   z_6 \right)\delta+2 z_7  \delta
   '  \\\nonumber \{z_2(x),z_2(y)\}&=&
   \delta ' ,~\{z_2(x),z_4(y)\}\,=\, -\sqrt{2} z_7  \delta  ,~\{z_2(x),z_5(y)\}\,=\, \sqrt{2} z_6  \delta
   \\\nonumber \{z_2(x),z_6(y)\}&=& \sqrt{2} z_5  \delta  ,~~\{z_2(x),z_7(y)\}\,=\, -\sqrt{2} z_4  \delta
    \\\nonumber \{z_3(x),z_3(y)\}&=&
   \left(z_1' -3
   z_5  z_5' +3
   z_6
   z_6' \right)\delta  + \left
   (-3 z_5 ^2+3 z_6 ^2+2
   z_1 \right) \delta
   ' -\frac{1}{2}\delta ^{(3)}   \\\nonumber
    \{z_3(x),z_4(y)\}&=& (-2 z_6 z_2'-z_2
   z_6'+\frac{1}{\sqrt{2}}z_5''+\sqrt{2} z_5^3+\sqrt{2} z_2^2 z_5-\sqrt{2} z_6^2
   z_5-2 \sqrt{2} z_1 z_5)\delta \\\nonumber & & +(\frac{3}{\sqrt{2}} z_5'-3 z_2 z_6)\delta'+\frac{3}{\sqrt{2}} z_5\delta''
   \\\nonumber \{z_3(x),z_5(y)\}&=& -\sqrt{2} z_4  \delta  ,~  \{z_3 (x),z_6(y)\}\,=\, \sqrt{2} z_7  \delta
     \\\nonumber
     \{z_3(x) ,z_7(y)\}&=& (2 z_5z_2'+z_2
   z_5'-\frac{1}{\sqrt{2}}z_6''+\sqrt{2} z_6^3-\sqrt{2} z_2^2 z_6-\sqrt{2} z_5^2
   z_6+2 \sqrt{2} z_1 z_6)\delta \\\nonumber & & + (3 z_2 z_5-\frac{3}{\sqrt{2}} z_6')\delta'-\frac{3}{\sqrt{2}} z_6\delta''
      \\\nonumber \{z_4(x),z_4(y)\}&=&
   \left(- z_1'+3
   z_2 z_2'+3
   z_5
   z_5'\right) \delta +\left
   (3z_2^2+3z_5^2-2 z_1\right) \delta
   '+\frac{1}{2} \delta ^{(3)}
        \\\nonumber \{z_4(x),z_7(y)\}&=&  (-z_6 z_5'-2
   z_5
   z_6'+\frac{1}{\sqrt{2}}
   z_2''+\sqrt{
   2} z_2^3+\sqrt{2} z_5^2
   z_2-\sqrt{2} z_6^2
   z_2-2 \sqrt{2} z_1 z_2)\delta
   \\\nonumber
  & &+ (\frac{3}{\sqrt{2}} z_2'-3 z_5 z_6)\delta'+ \frac{3 }{\sqrt{2}}z_2 \delta''
   \\\nonumber  \{z_4(x),z_5(y)\}&=& -\sqrt{2} z_3  \delta,~ \{z_5(x),z_5(y)\}= \delta ',~
         \{z_5(x),z_6(y)\}= -\sqrt{2} z_2  \delta
       \\\nonumber \{z_6(x),z_6(y)\}&=& -\delta ',~
      \{z_6(x),z_7(y)\}= \sqrt{2} z_3  \delta
       \\\nonumber \{z_7(x) ,z_7(y)\}&=& (z_1'-3 z_2 z_2'+3 z_6
   z_6')\delta+ (-3 z_2^2+3 z_6^2+2 z_1)\delta'-\frac{1}{2}\delta^{(3)},
 \end{eqnarray}
and  $\mathbb B_1^\Q=\Lie_{\partial_{z_3}} \mathbb B_2^\Q$. The leading terms of the local Poisson brackets reads 

{\small \[F^{ij}_2=\sqrt{2}\left(
\begin{array}{ccccccc}
 0 & 0 & -z_4 z_5-z_6 z_7 & z_2 z_7-z_3 z_5 & 0 & 0 & z_2 z_4+z_3
   z_6 \\
\ * & 0 & 0 & -z_7 & z_6 & z_5 & -z_4 \\
\ * & * & 0 & z_5^3+z_2^2 z_5-z_6^2 z_5-2 z_1 z_5 &
   -z_4 & z_7 & z_6^3-z_2^2 z_6-z_5^2 z_6+2 z_1 z_6 \\
\ * & * & * & 0
   & -z_3 & 0 & z_2^3+z_5^2 z_2-z_6^2 z_2-2 z_1 z_2 \\
\ * & * & * & * & 0 & -z_2 & 0 \\
\ * & * & * & * & * & 0 & z_3 \\
\ * & * & * & * & * & * & 0
\end{array}
\right) \]}

\[\Omega^{ij}_2=\left(
\begin{array}{ccccccc}
 2 z_1 & 0 & 2 z_3 & 2 z_4 & 0 & 0 & 2 z_7 \\
 0 & 1 & 0 & 0 & 0 & 0 & 0 \\
 2 z_3 & 0 & -3 z_5^2+3 z_6^2+2 z_1 & -3 z_2 z_6 & 0 & 0 & 3 z_2
   z_5 \\
 2 z_4 & 0 & -3 z_2 z_6 & 3 z_2^2+3 z_5^2-2 z_1 & 0 & 0 & -3 z_5
   z_6 \\
 0 & 0 & 0 & 0 & 1 & 0 & 0 \\
 0 & 0 & 0 & 0 & 0 & -1 & 0 \\
 2 z_7 & 0 & 3 z_2 z_5 & -3 z_5 z_6 & 0 & 0 & -3 z_2^2+3 z_6^2+2
   z_1
\end{array}
\right)\]
and \[ S_{1;2}^{ij}=\left(
\begin{array}{ccccccc}
 0 & 0 & 0 & 0 & 0 & 0 & 0 \\
 0 & 0 & 0 & 0 & 0 & 0 & 0 \\
 0 & 0 & 0 & \frac{3}{\sqrt{2}} z_5 & 0 & 0 & -\frac{3
   }{\sqrt{2}} z_6\\
 0 & 0 & -\frac{3 }{\sqrt{2}}z_5 & 0 & 0 & 0 & \frac{3
   }{\sqrt{2}}z_2 \\
 0 & 0 & 0 & 0 & 0 & 0 & 0 \\
 0 & 0 & 0 & 0 & 0 & 0 & 0 \\
 0 & 0 & \frac{3}{\sqrt{2}} z_6 & -\frac{3}{\sqrt{2}} z_2 & 0 &
   0 & 0
\end{array}
\right)\]

The space of common equilibrium points $N$ is defined by
\be \label{case22}
N:=\{q\in Q: \partial_{z_\alpha}\overline P_j^0=0, ~\alpha=4,5,6,7, j=1,2,3\}=\{q\in Q: z_4=z_5=z_5=z_7=0\}
\ee
Consider the reduced  bihamiltonian structure on loop space $\N$ which admits a dispersionless limit.  Thus the matrix $\Omega^{ab}_2$ is found by taking the upper left minor of $\Omega^{ij}_2$ under  the restriction  \eqref{case22}.  Also in this case,  the matrices of the leading terms $\Omega^{ab}_2$ and $\Omega_1^{ab}=\partial_{z_3}\Omega^{ab}_2$ do not define a flat pencil of metrics   since $\det \Omega_1^{ab}=0$. However, the matrix $S_{2;1}^{ab}=0$.

As explained above, we seek a vector field $e$ such that $\Lie^{2}_{e}  \Omega^{ab}_2=0$ and  $ \Lie_{e}  \Omega^{ab}_2$ define a flat contravariant metric. By  ad-hoc trials, we find
\[e= -\partial_{z_1}+  z_2^{-1} \partial_{z_2}\]
satisfies the required condition. Moreover,  $\Lie_{e}^2 \mathbb B_2^\N=0$. Hence, we get a new local bihamiltonian structure on $\N$ formed by $\mathbb B_2^\N$ and $\overline \mathbb B_1^\N:=\Lie_{e} \mathbb B_2^\N$ admitting a dispersionless limit. Moreover, the leading terms define a regular quasihomogeneous flat pencil of metrics formed by $\Omega^{ab}_2$ and  $\overline\Omega^{ab}_1:=\Lie_e\Omega^{ab}_2$ of degree 0.

In the  coordinates
\[t_1=z_1+\frac{1}{2}z_2^2,~t_2=\frac{z_2^2}{2},~t_3=z_3. \]

We get 

\[\Omega_2^{ab}(t)=\left(
\begin{array}{ccc}
 2 t_1 & 2 t_2 & 2 t_3 \\
 2 t_2 & 2 t_2 & 0 \\
 2 t_3 & 0 & 2 t_1-2 t_2
\end{array}
\right), ~S_{1;2}^{ab}=0 ,~\widetilde\Omega_1^{ab}(t)=\left(
\begin{array}{ccc}
 0 & 2 & 0 \\
 2 & 2 & 0 \\
 0 & 0 & -2
\end{array}
\right)\]
and $e=\partial_{t_2}$. The resulting Dubrovin-Frobenius manifold structure is trivial, i.e.,   Frobenius algebras on the tangent space are independent of the point. It has the potential
\[ \mathbb F=\frac{t_1^3}{12}-\frac{1}{4} t_2 t_1^2+\frac{1}{4} t_2^2
   t_1+\frac{1}{4} t_3^2 t_1-\frac{1}{4} t_2 t_3^2.\]
The  Euler vector field is $E=\sum t_i\partial_{t_i}$ satisfying $E\mathbb F=3 \mathbb F$. To find the central invariants, we have
\[S_{2;2}^{ab}(t)=\left(
\begin{array}{ccc}
 -\frac{1}{2} & 0 & 0 \\
 0 & 0 & 0 \\
 0 & 0 & -\frac{1}{2}
\end{array}
\right).\]
Then using equations \eqref{cent form}, the central invariants of the associated bihamiltonian structure are $0$, $-\frac{1}{48}$ and $-\frac{1}{48}$. In particular, the bihamiltonian structure on $\N$ is not of topological type. We have no explanation for the occurrence of uneven constant central invariants in this case.

\subsection{Conclusion}

We give examples of bihamiltonian structure of topological type associated to nilpotent elements of semisimple  type beginning form the associated classical $W$-algebras.  For subregular nilpotent element in $sl_3$ and $sl_4$ the corresponding Dubrovin-Frobenius manifold is logarithmic. Moreover, it was proved in \cite{LZ} that the bihamiltonian structures are associated to  constrained KP hierarchies. It seems this overlap is  true for other subregular nilpotent elements in $sl_n$, $n>2$. In future work, using classical $W$-algebras, we will give uniform construction for  bihamiltonian structures of topological type for  logarithmic Dubrovin-Frobenius manifolds  and analyze their connection to bihamiltonian structures associated to the constrained KP hierarchies.

\vspace{0.1cm}

\noindent{\bf Acknowledgments.}
 The author thanks Marco Pedroni for his time to read the first draft of this article and for his useful suggestion to include the example of bihamiltonian structure associated to the fractional KDV hierarchy and brings to my attention the paper \cite{LZ}. The algebraic Dubrovin-Frobenius structure  given there was first computed by Yougin Zhang in his visit to Sultan Qaboos University in 2017. The author would like to thank the Isaac Newton Institute for Mathematical Sciences for support and hospitality during the programme Dispersive Hydrodynamics when work on this paper was undertaken (EPSRC Grant Number EP/R014604/1).
\vspace{0.1cm}
 
\noindent{\bf Funding} This work was partially funded by the internal grant of Sultan Qaboos University \\ (IG/SCI/DOMS/19/08).

\vspace{0.1cm}

\noindent{\bf Data Availability} Non applicable. 
\section*{Declarations}

\noindent{\bf Conflict of interest} The authors have no relevant financial or non-financial interests to disclose.

\noindent Yassir Dinar

\noindent Department of Mathematics

\noindent College of Science

\noindent Sultan Qaboos University, Oman

\noindent dinar@squ.edu.om.

\begin{thebibliography}{99}
\bibitem{bolv1}
Bolsinov, A.V., Oshemkov, A.A.; Bi-Hamiltonian structures and singularities of integrable systems. Regul. Chaot. Dyn. 14, 431–454 (2009).

\bibitem{gDSh2}  Burroughs, N., de Groot, M., Hollowood, T. and Miramontes, J.;
 Generalized Drinfeld-Sokolov hierarchies II: the Hamiltonian
 structures, Comm. Math. Phys.153, 187 (1993).
\bibitem{COLMC} Collingwood, D. H., McGovern, W. M.; Nilpotent orbits in semisimple Lie algebras. Van Nostrand Reinhold Mathematics Series. ISBN: 0-534-18834-6 (1993).
\bibitem{DelFeher} Delduc, F.; Feher, L., Regular conjugacy classes in the Weyl group and integrable hierarchies.  J. Phys. A  28,  no. 20, 5843--5882  (1995).

\bibitem{mypaper6} Dinar, Yassir; Algebraic classical W-algebras and Frobenius manifolds. Lett Math Phys 111, 115 (2021).
\bibitem{mypaper4} Dinar, Yassir; $W$-algebras and the equivalence of bihamiltonian, Drinfeld-Sokolov and Dirac reductions. J. Geom. Phys. 84, 30-42 (2014).

\bibitem{mypaper5} Dinar, Y.; On integrability of transverse Lie-Poisson structure to nilpotent elements. Journal of Geometry and Physics, Volume 155, 103690, ISSN 0393-0440 (2020).

\bibitem{wdvv} Dijkgraaf, R., Verlinde, H. and Verlinde, E.; Topological strings in $d  1$.
Nucl. Phys. B352, 59 (1991).

\bibitem{DS} Drinfeld, V. G., Sokolov, V. V.; Lie algebras and equations of Korteweg-de Vries type. (Russian)  Current problems in mathematics, Vol. 24,  81--180, Itogi Nauki i Tekhniki, Akad. Nauk SSSR, Vsesoyuz. Inst. Nauchn. i Tekhn. Inform., Moscow, (1984).

 \bibitem{DuRev} Dubrovin, Boris; Geometry of $2$D topological field theories.  Integrable systems and quantum groups (Montecatini Terme, 1993),  120--348, Lecture Notes in Math., 1620, Springer, Berlin, (1996).


   \bibitem{DFP} Dubrovin, Boris; Flat pencils of metrics and Frobenius manifolds.  Integrable systems and algebraic geometry (Kobe/Kyoto, 1997),  47--72, World Sci. Publ. (1998).

   \bibitem{DZ} Dubrovin, B. , Zhang, Y.; Normal forms of hierarchies of integrable PDEs, Frobenius
manifolds and Gromov-Witten invariants, www.arxiv.org
math/0108160.
\bibitem{DLZ} Dubrovin, B. , Liu, Si-Qi,  Zhang, Y.; Frobenius manifolds and central invariants for the Drinfeld–Sokolov bihamiltonian structures, Advances in Mathematics, Volume 219, Issue 3,  780-837 (2008).
\bibitem{DN} Dubrovin, B. A., Novikov, S. P.; Poisson brackets of hydrodynamic type. (Russian)  Dokl. Akad. Nauk SSSR  279,  no. 2, 294--297  (1984).
\bibitem{DLiZ} Dubrovin, B., Liu, S.-Q., Zhang, Y.; On Hamiltonian perturbations of hyperbolic systems of conservation laws. I. Quasi-triviality of bi-Hamiltonian perturbations. Comm. Pure Appl. Math. 59, 559--615 (2006).
\bibitem{DF}  Dubrovin, B.; On universality of critical behaviour in Hamiltonian PDEs, Amer. Math. Soc. Transl. 224, pp 59–109 (2008). 
  \bibitem{Elash} Elashvili, A. G.; Kac, V. G.; Vinberg, E. B.; Cyclic elements in semisimple Lie algebras. Transform. Groups 18, no. 1,97-130 (2013).
  \bibitem{FalLor} Falqui, G. , Lorenzoni, P.; Exact Poisson pencils, $\tau$-structures and topological hierarchies. Phys. D 241 , no. 23-24, 2178–2187 (2012).  

  \bibitem{fehercomp} Feher, L., O'Raifeartaigh, L., Ruelle, P., Tsutsui, I.; On the completeness of the set of classical $ W$-algebras obtained from DS reductions.  Comm. Math. Phys.  162 ,  no. 2, 399--431 (1994).
\bibitem{LZ1}  Liu, S.-Q., Zhang, Y.: Deformations of semisimple bihamiltonian structures of hydrodynamic type. J. Geom. Phys. 54, 427--453 (2005).
\bibitem{LZ} Liu, S., Zhang, Y., Zhou, Xu; Central invariants of the constrained KP hierarchies. J. Geom. Phys. 97, 177-189 (2015).

\bibitem{Lo} Lorenzoni, P.; Deformations of bi-Hamiltonian structures of hydrodynamic type. J. Geom. Phys. 44, 331--375 (2002).

\bibitem{Ser} Sergyeyev, A. A; Simple Way of Making a Hamiltonian System Into a Bi-Hamiltonian One. Acta Applicandae Mathematicae 83, 183–197 (2004). 

\end{thebibliography}
\end{document}